\documentclass[11pt]{article}

\usepackage{amsmath,amsthm,amsfonts,bbold}

\usepackage{amssymb}

\usepackage{epsfig}

\usepackage{rotating}

\usepackage{subfigure}

\begin{document}

\title{Conditional path sampling for stochastic differential equations through drift relaxation}
\author{Panos Stinis \\ 
Department of Mathematics \\
University of Minnesota \\
    Minneapolis, MN 55455} 
\date {}

\maketitle

\begin{abstract}
We present an algorithm for the efficient sampling of conditional paths of stochastic differential equations (SDEs). While unconditional path sampling of SDEs is straightforward, albeit expensive for high dimensional systems of SDEs, conditional path sampling can be difficult even for low dimensional systems. This is because we need to produce sample paths of the SDE which respect both the dynamics of the SDE and the initial and endpoint conditions. The dynamics of a  SDE are governed by the deterministic term (drift) and the stochastic term (noise). Instead of producing conditional paths directly from the original SDE, one can consider a sequence of SDEs with modified drifts. The modified drifts should be chosen so that it is easier to produce sample paths which satisfy the initial and endpoint conditions.  Also, the sequence of modified drifts converges to the drift of the original SDE. We construct a simple Markov Chain Monte Carlo (MCMC) algorithm which samples, in sequence, conditional paths from the modified SDEs, by taking the last sampled path at each level of the sequence as an initial condition for the sampling at the next level in the sequence. The algorithm can be thought of as a stochastic analog of deterministic homotopy methods for solving nonlinear algebraic equations or as a SDE generalization of simulated annealing. The algorithm is particularly suited for filtering/smoothing applications. We show how it can be used to improve the performance of particle filters. Numerical results for filtering of a stochastic differential equation are included. 
\end{abstract}


\section*{Introduction}
The study of systems arising in different areas, from signal processing and chemical kinetics to econometrics and finance (see e.g. \cite{bolhuis,stuart}) often requires the sampling of paths of stochastic differential equations (SDEs) subject to initial and endpoint conditions. While unconditional path sampling of SDEs is straightforward, albeit expensive for high dimensional systems of SDEs, conditional path sampling can be difficult even for low dimensional systems. This is because we need to produce sample paths of the SDE which respect both the dynamics of the SDE and the initial and endpoint conditions. An analogous situation arises in ordinary differential equations, where it can be considerably more difficult to create solutions to boundary value problems than it is to construct solutions to initial value problems (see e.g. Ch. 8 in \cite{deuflhard}). The problem of conditional path sampling of SDEs has been a subject of active research in recent years and some very interesting approaches have already been developed (see e.g. \cite{C08,stuart,W07}). 

The dynamics of a SDE are governed by the deterministic term (drift) and the stochastic term (noise). Instead of producing conditional paths directly from the original SDE, one can consider a sequence of SDEs with modified drifts. The modified drifts should be chosen so that it is easier to produce sample paths which satisfy the initial and endpoint conditions.  Also, the sequence of modified drifts converges to the drift of the original SDE. We construct a simple Markov Chain Monte Carlo (MCMC) algorithm which samples, in sequence, conditional paths from the modified SDEs, by taking the last sampled path at each level of the sequence as an initial condition for the sampling at the next level in the sequence. 

We have used the drift relaxation algorithm to modify a popular filtering method called particle filter \cite{doucet}. A particle filter is a sequential importance sampling algorithm which is based on the recursive (online) Bayesian updating of the values of samples (called particles) to incorporate information from noisy observations of the state of a dynamic model. While the particle filter is a very versatile method it may require a very large number of samples to approximate accurately the conditional density of the state of the model. This has led to considerable research (see e.g. \cite{gilks,W09}) into how one can modify a particle filter to make it more efficient (see also \cite{CT09} for a different approach to particle filtering). As an application of the drift relaxation algorithm we show in Section \ref{particle_filtering} how it can be used to construct a more efficient particle filter.     

The paper is organized as follows. Section \ref{drift-relaxation} presents the drift relaxation algorithm for an SDE conditional path sampling problem. Section \ref{particle_filtering} shows how to use the algorithm to modify a particle filter. Section \ref{numerical} contains numerical results for the application of the modified particle filter to the standard example of filtering a diffusion in a double-well potential (more elaborate examples will be presented in \cite{MS10}). Finally, Section \ref{discussion} discusses the results as well as current and future work.

\section{Conditional path sampling and drift relaxation}\label{drift-relaxation}

Suppose that we are given a system of stochastic differential equations (SDEs)
\begin{equation}\label{sdes}
dX_t= a(X_t)dt + \sigma(X_t) dB_t,
\end{equation}
Suppose also that we want to construct, in the time interval $[0,T],$ sample paths from \eqref{sdes} such that the endpoints are distributed according to the densities $h(X_0)$ and $g(X_T)$ respectively. Equation \eqref{sdes} can be discretized in the interval $[0,T]$ by some numerical approximation scheme \cite{kloeden}. Suppose that we have discretized the interval $[0,T]$ using a stepsize $\Delta t= T/I.$ To construct conditional paths of \eqref{sdes} we can start by sampling initial conditions from the density $h.$ For each initial condition $X_0$ we can then produce the desired conditional paths by sampling the density 
\begin{equation}\label{path-density}
\prod_{i=1}^{I}p(X_{T_i}|X_{T_{i-1}})g(X_T),
\end{equation}
where $p(X_{T_i}|X_{T_{i-1}})$ is the transition probability from $X_{T_{i-1}}$ at time $T_{i-1}$ to the point $X_{T_i}$ at time $T_i.$ The density \eqref{path-density} can be sampled using MCMC assuming that the transition densities $p(X_{T_i}|X_{T_{i-1}})$ can be evaluated. However, the major issue with the MCMC sampling is whether it can be performed efficiently (see e.g. \cite{W07,CT09}). Instead of MCMC sampling directly from the density \eqref{path-density} i.e., starting from an arbitrary initial path and modifying it to become a path corresponding to \eqref{path-density}, we can aid the MCMC sampling process by providing the MCMC sampler of the density \eqref{path-density} with a better initial condition.

To this end, consider an SDE system with modified drift 
\begin{equation}\label{sdes-simple}
dY_t= b(Y_t)dt + \sigma(Y_t) dB_t,
\end{equation}
where $b(Y_t)$ can be suitably chosen to facilitate the conditional path sampling problem.

Also, consider the collection of $L+1$ modified SDE systems
\begin{equation*}
dY^l_t= (1-\epsilon_l) b(Y^l_t)dt + \epsilon_l a(Y^l_t) dt + \sigma(Y^l_t) dB_t,
\end{equation*}
where $\epsilon_l  \in [0,1], \; l=0,\ldots,L,$ with $\epsilon_l < \epsilon_{l+1},$ $\epsilon_0=0$ and $\epsilon_{L}=1.$ Instead of sampling directly a (conditional) path from the SDE \eqref{sdes}, one can sample a path from the modified SDE \eqref{sdes-simple} and gradually morph the path 
into a path of \eqref{sdes}.

\vskip 14pt {\bf Drift relaxation algorithm:}
\begin{itemize}
\item
($l=0$) Begin with a sample path from the modified SDE \eqref{sdes-simple}.
\item
Sample through MCMC the density $\prod_{i=1}^{I}p(Y^0_{T_i}|Y^0_{T_{i-1}})g(Y^0_T).$
\item
For $l=1,...,L$ take the last sample path from the ($l-1$)st SDE and use it as in initial condition for MCMC sampling the density 
$$\prod_{i=1}^{I}p(Y^l_{T_i}|Y^l_{T_{i-1}})g(Y^l_T) $$
at the $l$th level. 
\item
Keep the last sample path at the $L$th level.
\end{itemize}
The drift relaxation algorithm is similar to Simulated Annealing (SA) used in equilibrium statistical mechanics \cite{liu}. However, instead of modifying a temperature as in SA, here we modify the drift of the system.

\section{Application to particle filtering}\label{particle_filtering}

We show in this section how the drift relaxation algorithm can be applied to particle filtering with the aim of bringing the samples closer to the observations.  

\subsection{Generic particle filter}\label{generic}

Suppose that we are given an SDE system and that we also have
access to noisy observations $Z_{T_1},\ldots,Z_{T_K}$ of the state
of the system at specified instants $T_1,\ldots,T_K.$ The
observations are functions of the state of the system, say given
by $Z_{T_k}=G(X_{T_k},\xi_k),$ where $\xi_k, k=1,\ldots,K$ are
mutually independent random variables. For simplicity, let us
assume that the distribution of the observations admits a density
$g(X_{T_k},Z_{T_k}),$ i.e., $p(Z_{T_k}|X_{T_k} ) \propto
g(X_{T_k},Z_{T_k}).$

The filtering problem consists of computing estimates of the
conditional expectation $E[f(X_{T_k})| \{Z_{T_j}\}^{k}_{j=1}],$
i.e., the conditional expectation of the state of the system given
the (noisy) observations. Equivalently, we are looking to compute
the conditional density of the state of the system given the
observations $p(X_{T_k}|\{Z_{T_j}\}^{k}_{j=1}).$ There are several
ways to compute this conditional density and the associated
conditional expectation but for practical applications they are
rather expensive.

Particle filters fall in the category of importance sampling
methods. Because computing averages with respect to the
conditional density involves the sampling of the conditional
density which can be difficult, importance sampling methods
proceed by sampling a reference density
$q(X_{T_k}|\{Z_{T_j}\}^{k}_{j=1})$ which can be easily sampled and
then compute the weighted sample mean
$$E[f(X_{T_k})| \{Z_{T_j}\}^{k}_{j=1}] \approx \frac{1}{N} \sum_{n=1}^N
f(X^n_{T_k})\frac{p(X^n_{T_k}|\{Z_{T_j}\}^{k}_{j=1})}{q(X^n_{T_k}|\{Z_{T_j}\}^{k}_{j=1})}$$
or the related estimate
\begin{equation}\label{importance}
E[f(X_{T_k})| \{Z_{T_j}\}^{k}_{j=1}] \approx \frac{\sum_{n=1}^N
f(X^n_{T_k})
\frac{p(X^n_{T_k}|\{Z_{T_j}\}^{k}_{j=1})}{q(X^n_{T_k}|\{Z_{T_j}\}^{k}_{j=1})}}{\sum_{n=1}^N
\frac{p(X^n_{T_k}|\{Z_{T_j}\}^{k}_{j=1})}{q(X^n_{T_k}|\{Z_{T_j}\}^{k}_{j=1})}},
\end{equation}
where $N$ has been replaced by the approximation
$$N \approx \sum_{n=1}^N \frac{p(X^n_{T_k}|\{Z_{T_j}\}^{k}_{j=1})}{q(X^n_{T_k}|\{Z_{T_j}\}^{k}_{j=1})}.$$
Particle filtering is a recursive implementation of the importance
sampling approach. It is based on the recursion
\begin{align}
p(X_{T_k}|\{Z_{T_j}\}^{k}_{j=1}) & \propto g(X_{T_k},Z_{T_k}) p(X_{T_k}|\{Z_{T_j}\}^{k-1}_{j=1}),
\label{correct}\\  
\text{where} \;\; p(X_{T_k}|\{Z_{T_j}\}^{k-1}_{j=1}) & =  \int p(X_{T_k}|
X_{T_{k-1}})p(X_{T_{k-1}}|\{Z_{T_j}\}^{k-1}_{j=1}) dX_{T_{k-1}}.
\label{update}
\end{align}
If we set
$$q(X_{T_k}|\{Z_{T_j}\}^{k}_{j=1})=p(X_{T_k}|\{Z_{T_j}\}^{k-1}_{j=1}),$$
then from \eqref{correct} we get
$$\frac{p(X_{T_k}|\{Z_{T_j}\}^{k}_{j=1})}{q(X_{T_k}|\{Z_{T_j}\}^{k}_{j=1})} \propto g(X_{T_k},Z_{T_k}).$$
The approximation in expression \eqref{importance} becomes
\begin{equation}\label{particle}
E[f(X_{T_i})| \{Z_{T_j}\}^{k}_{j=1}] \approx \frac{\sum_{n=1}^N
f(X^n_{T_k})g(X^n_{T_k},Z_{T_k})}{\sum_{n=1}^N
g(X^n_{T_k},Z_{T_k})}
\end{equation}
From \eqref{particle} we see that if we can construct samples from
the predictive distribution $p(X_{T_k}|\{Z_{T_j}\}^{k-1}_{j=1})$
then we can define the (normalized) weights $W^n_{T_k}=
\frac{g(X^n_{T_k},Z_{T_k})}{\sum_{n=1}^N g(X^n_{T_k},Z_{T_k})},$
use them to weigh the samples and the weighted samples will be
distributed according to the posterior distribution
$p(X_{T_k}|\{Z_{T_j}\}^{k}_{j=1}).$

In many applications, most samples will have a negligible weight
with respect to the observation, so carrying them along does not
contribute significantly to the conditional expectation estimate
(this is the problem of degeneracy \cite{liu}). To create larger
diversity one can resample the weights to create more copies of
the samples with significant weights. The particle filter with
resampling is summarized in the following algorithm due to Gordon
{\it et al.} \cite{gordon}.

\vskip14pt
{\bf Particle filter}
\begin{enumerate}
\item Begin with $N$ unweighted samples $X^n_{T_{k-1}}$ from
$p(X_{T_{k-1}}|\{Z_{T_j}\}^{k-1}_{j=1}).$ \item {\bf Prediction}:
Generate $N$ samples $X'^n_{T_k}$ from $ p(X_{T_k}| X_{T_{k-1}}).$
\item {\bf Update}: Evaluate the weights $$W^n_{T_k}=
\frac{g(X'^n_{T_k},Z_{T_k})}{\sum_{n=1}^N
g(X'^n_{T_k},Z_{T_k})}.$$ 
\item {\bf Resampling}: Generate $N$
independent uniform random variables $\{\theta^n\}_{n=1}^N$ in
$(0,1).$ For $n=1,\ldots,N$ let $X^n_{T_k}=X'^j_{T_k} $where $$
\sum_{l=1}^{j-1}W^l_{T_k} \leq \theta^j < \sum_{l=1}^{j}W^l_{T_k}
$$ 
where $j$ can range from $1$ to $N.$
\item Set $k=k+1$ and proceed to Step 1.
\end{enumerate}

The particle filter algorithm is easy to implement and adapt for
different problems since the only part of the algorithm that
depends on the specific dynamics of the problem is the prediction
step. This has led to the particle filter algorithm's increased
popularity \cite{doucet}. However, even with the resampling step,
the particle filter can still need a lot of samples in order to
describe accurately the conditional density
$p(X_{T_k}|\{Z_{T_j}\}^{k}_{j=1}).$ Snyder {\it et al.}
\cite{snyder} have shown how the particle filter can fail in
simple high dimensional problems because one sample dominates the
weight distribution. The rest of the samples are not in
statistically significant regions. Even worse, as we will show in
the numerical results section, there are simple examples where not
even one sample is in a statistically significant region. In the
next subsection we present how drift relaxation can be used to
push samples closer to statistically significant regions.

\subsection{Particle filter with MCMC step}\label{MCMC_step}
Several authors (see e.g. \cite{gilks,W09}) have suggested the use
of a MCMC step after the resampling step (Step 4) in order to move
samples away from statistically insignificant regions. There are
many possible ways to append an MCMC step after the resampling
step in order to achieve that objective. The important point is
that the MCMC step must preserve the conditional density
$p(X_{T_k}|\{Z_{T_j}\}^{k}_{j=1}).$ In the current section we show
that the MCMC step constitutes a case of conditional path
sampling.

We begin by noting that one can use the resampling step (Step 4)
in the particle filter algorithm to create more copies not only of
the good samples according to the observation, but also of the
values (initial conditions) of the samples at the previous
observation. These values are the ones who have evolved into good
samples for the current observation (see more details in
\cite{W09}). The motivation behind producing more copies of the
pairs of initial and final conditions is to use the good initial
conditions as starting points to produce statistically more
significant samples according to the current observation. This
process can be accomplished in two steps. First, Step 4 of the
particle filter algorithm is replaced by

\vskip14pt {\bf Resampling}: Generate $N$ independent uniform
random variables $\{\theta^n\}_{n=1}^N$ in $(0,1).$ For
$n=1,\ldots,N$ let
$(X^n_{T_{k-1}},X^n_{T_k})=(X'^j_{T_{k-1}},X'^j_{T_k}) $where $$
\sum_{l=1}^{j-1}W^l_{T_k} \leq \theta^j < \sum_{l=1}^{j}W^l_{T_k}
$$
Also, through Bayes' rule \cite{W09} one can show that the posterior density $p(X_{T_k}|\{Z_{T_j}\}^{k}_{j=1})$ is 
preserved if one samples from the density $$g(X_{T_k},Z_{T_k}) p(X_{T_k}|X_{T_{k-1}})$$
where $X_{T_{k-1}}$ are given by the modified resampling step. This is a problem of conditional path sampling for (continuous-time or discrete) stochastic systems. The important issue is to perform the necessary sampling efficiently \cite{CT09,W09}. We propose to do that here using drift relaxation (see Section \ref{drift-relaxation}). The particle filter with MCMC step algorithm is given by

\vskip14pt
{\bf Particle filter with MCMC step}
\begin{enumerate}
\item Begin with $N$ unweighted samples $X^n_{T_{k-1}}$ from
$p(X_{T_{k-1}}|\{Z_{T_j}\}^{k-1}_{j=1}).$ 
\item {\bf Prediction}:
Generate $N$ samples $X'^n_{T_k}$ from $ p(X_{T_k}| X_{T_{k-1}}).$
\item {\bf Update}: Evaluate the weights $$W^n_{T_k}=
\frac{g(X'^n_{T_k},Z_{T_k})}{\sum_{n=1}^N
g(X'^n_{T_k},Z_{T_k})}.$$ 
\item {\bf Resampling}: Generate $N$
independent uniform random variables $\{\theta^n\}_{n=1}^N$ in
$(0,1).$ For $n=1,\ldots,N$ let
$(X^n_{T_{k-1}},X^n_{T_k})=(X'^j_{T_{k-1}},X'^j_{T_k}) $ where $$
\sum_{l=1}^{j-1}W^l_{T_k} \leq \theta^j < \sum_{l=1}^{j}W^l_{T_k}
$$
where $j$ can range from $1$ to $N.$
 \item {\bf MCMC step}: For $n=1,\ldots,N$ choose a modified
drift (possibly different for each $n$). Construct a path for the
SDE with the modified drift starting from $X^n_{T_{k-1}}.$ Construct through drift relaxation a Markov chain
for $Y^{n}_{T_k}$ with stationary distribution
$$ g(Y^n,Z_{T_k}) p(Y^n | X^n_{T_{k-1}}) $$
\item
Set $X^n_{T_k}=Y^{n}_{T_k}.$
\item
Set $k=k+1$ and proceed to Step 1.
\end{enumerate}
Since the samples $X^n_{T_k}=Y^{n}_{T_k}$ are constructed
by starting from different sample paths, they are independent.
Also, note that the samples $X^n_{T_k}$ are unweighted. However, we
can still measure how well these samples approximate the posterior
density by comparing the effective sample sizes of the particle
filter with and without the MCMC step. For a collection of $N$
samples the effective sample size $ess(T_k)$ is defined by
$$ess(T_k) = \frac{N}{1+C_k^2}$$
where
\begin{equation*}
C_k =  \frac{1}{W_k}
\sqrt{\frac{1}{N}\sum_{n=1}^N (g(X^n_{T_k},Z_{T_k}) -W_k)^2} \;\; \text{and} \;\; W_k = \frac{1}{N} \sum_{n=1}^N g(X^n_{T_k},Z_{T_k}).
\end{equation*}
The effective sample size can be interpreted as that the $N$ weighted samples are worth of $ess(T_k) = \frac{N}{1+C_k^2}$ i.i.d. samples drawn from the target density, which in our case is the posterior density. By definition, $ess(T_k) \le N.$ If the samples have uniform weights, then $ess(T_k)=N.$ On the other hand, if all samples but one have zero weights, then $ess(T_k)=1.$

\section{Numerical results}\label{numerical}
We present numerical results of the particle filter algorithm with MCMC step for the standard problem of diffusion in a double-well potential (more elaborate applications of the method will be presented elsewhere \cite{MS10}). Our objective here is to show how the generic particle filter's performance can be significantly improved by incorporating the MCMC step via drift relaxation.

The problem of diffusion in a double well potential is described by the scalar SDE
\begin{equation}\label{double}
dX_t=-4X_t(X_t^2-1) + \frac{1}{2} dB_t.
\end{equation}
The deterministic part (drift) describes a gradient flow for the potential $U(x)=x^4-2x^2$ which has two minima, at $x=\pm1.$ In the notation of Section \ref{particle_filtering} we have $a(X_t)=-4X_t(X_t^2-1)$ and $\sigma(X_t)=\frac{1}{2}.$ If the stochastic term is zero the solution wanders around one of the minima depending on the value of the initial condition. A weak stochastic term leads to rare transitions between the minima of the potential. We have chosen the coefficient $\frac{1}{2}$ to make the stochastic term rather weak. This is done because we plan to enforce the observations to alternate among the minima, and thus check if the particle filter can track these transitions.

The SDE \eqref{double} is discretized by the Euler-Maruyama \cite{kloeden} scheme with step size $\Delta t=10^{-2}$ which is small enough to guarantee stability of the scheme. The initial condition is set to $-1$ and there is a total of $10$ observations at $T_k=k, k=1,\ldots,10.$ The observations are given by $Z_{T_k}=X_{T_k} + \xi_k,$ where $\xi_k \sim N(0,0.01)$ for $k=1,\ldots,10.$ For this choice of observation noise, the observation density (also called likelihood) is given by 
\begin{equation}\label{double_observation}
g(X_{T_k},Z_{T_k}) \propto \exp \biggl[- \frac{( Z_{T_k}-X_{T_k})^2}{2*0.01}    \biggr] 
\end{equation}
The observations alternate between $1$ and $-1.$ In particular, for $k=1,\ldots,10$ we have $Z_{T_k}=-1$ if $k$ is odd and $Z_{T_k}=1$ if $k$ is even.

In order to apply the MCMC step with drift relaxation we need to define the modified drift $b(Y_t)$ for the process $Y_t$ given by 
\begin{equation}\label{double_simple}
dY_t=b(Y_t)+ \frac{1}{2}dB_t.
\end{equation}
The modified drift can be the same for all the samples or different for each sample. Since the difficulty in tracking the observations comes from the inability of the original SDE \eqref{double} to make frequent transitions between the two minima of the double well, an intuitively appealing choice for $b(Y_t)$ is $b(Y_t)=-\alpha4Y_t(Y_t^2-1),$ where $\alpha < 1.$ This drift corresponds to the potential $W(y)=\alpha(y^4-2y^2).$ The potential $W(y)$ has its minima also located at $y = \pm 1.$ However, the value of the potential at the minima is $-\alpha$ instead of $-1$ for the potential $U(x).$ This means that the wells corresponding to the minima of $W(y)$ are shallower than the wells corresponding to the minima of $U(x).$ This makes the transitions between the two wells for the process $Y_t$ more frequent than for the original process $X_t.$ For the numerical experiments we have chosen $\alpha=0.1.$

The sequence of modified SDEs for the drift relaxation algorithm with $L$ levels is given by 
\begin{equation}\label{sdes-relaxation}
dY^l_t= (1-\epsilon_l) b(Y^l_t)dt + \epsilon_l a(Y^l_t) dt + \frac{1}{2} dB_t,
\end{equation}
where $\epsilon_l  \in [0,1], \; l=0,\ldots,L,$ with $\epsilon_l < \epsilon_{l+1},$ $\epsilon_0=0$ and $\epsilon_{L}=1.$ For our numerical experiments we chose $L=10$ and $\epsilon_l = l/10.$

Recall that the density we want to sample during the MCMC step is given by $$g(X_{T_k},Z_{T_k}) p(X_{T_k}|X_{T_{k-1}}),$$ where $p(X_{T_k}|X_{T_{k-1}})$ is the transition probability between $X_{T_{k-1}}$ and $X_{T_k}.$ For many applications, sampling directly from $p(X_{T_k}|X_{T_{k-1}})$ may be impossible. Thus, one needs to resort to some numerical approximation scheme which approximates the path between $X_{T_{k-1}}$ and $X_{T_k}$ by a discretized path. However (see \cite{W09} for details), even the evaluation of the discretized path's density may not be efficient. Instead, by using the fact that each Brownian path in \eqref{double} gives rise to a unique path for $X_t$ \cite{oksendal}, we can replace the sampling of $g(X_{T_k},Z_{T_k}) p(X_{T_k}|X_{T_{k-1}})$ by sampling from the density
\begin{multline}\label{density}
\exp \biggl[- \frac{( Z_{T}-X^n_{T}(\{\Delta B^n_{i} \}_{i=0}^{I-1}))^2}{2*0.01}    \biggr]  \prod_{i=0}^{I-1} \exp \biggl[ - \frac{(\Delta B^n_{i})^2}{2*\Delta t}  \biggr] = \\
\exp \biggl[ - \biggl(  \frac{( Z_{T}-X^n_{T}(\{\Delta B^n_{i} \}_{i=0}^{I-1}))^2}{2*0.01}  +  \sum_{i=0}^{I-1} \frac{(\Delta B^n_{i})^2}{2*\Delta t}  \biggr) \biggr]
\end{multline}
where $\{\Delta B^n_{i} \}_{i=0}^{I-1}$ are the Brownian increments of the discretized path connecting  $X_{T_{k-1}}$ and $X_{T_k}.$ Also, note that the final point $X_{T_k}$ has now become a function of the entire Brownian path $\{\Delta B^n_{i} \}_{i=0}^{I-1}.$ For the numerical experiments we have chosen $\Delta t= \frac{T_k-T_{k-1}}{I}=10^{-2}$ which, since $T_k-T_{k-1}=1,$ gives $I=100.$

We use drift relaxation to produce samples from the density \eqref{density}. The Markov chain at each level of the drift relaxation algorithm is constructed using Hybrid Monte Carlo (HMC) \cite{liu}. At the $l$th level, we can discretize \eqref{sdes-relaxation}, say with the Euler-Maruyama scheme, and the points on the path will be given by 
$$Y^{l,n}_{i\Delta t}=Y^{l,n}_{(i-1)\Delta t} +  (1-\epsilon_l) b(Y^{l,n}_{(i-1)\Delta t})\Delta t + \epsilon_l a(Y^{l,n}_{(i-1)\Delta t})\Delta t + \frac{1}{2} \Delta B^{l,n}_{i-1},$$ 
for $i=1,\ldots,I.$ We can use more sophisticated schemes than the Euler-Maruyama scheme for the discretization of the simplified SDE \eqref{double_simple} at the cost of making the expression for the density more complicated.

We can define a potential $V_{\epsilon_l}(\{\Delta B^{l,n}_{i} \}_{i=0}^{I-1})$ for the variables $\{\Delta B^{l,n}_{i} \}_{i=0}^{I-1}.$ The potential is given by 
\begin{equation*}
V_{\epsilon_l} \bigl(\{\{\Delta B^{l,n}_{i} \}_{i=0}^{I-1} \bigr)= 
 \frac{( Z_{T}-Y^{l,n}_{I\Delta t}(\{\Delta B^{l,n}_{i} \}_{i=0}^{I-1}))^2}{2*0.01} 
+ \sum_{i=0}^{I-1} \frac{(\Delta B^{l,n}_{i})^2}{2*\Delta t}
\end{equation*}
and the density to be sampled can be written as $$\exp\biggl[-V_{\epsilon_l} \bigl(\{\Delta B^{l,n}_{i} \}_{i=0}^{I-1} \bigr)\biggr].$$ The subscript ${\epsilon_l}$ is to denote the dependence of the potential on the drift relaxation parameter ${\epsilon_l}.$ In HMC one considers the variables on which the potential depends as the position variables of a Hamiltonian system. In our case we have $I$ position variables so we can define a $I$-dimensional position vector $\{q_i\}_{i=1}^{I}.$ The next step is to augment the position variables vector by a vector of associated momenta $\{p_i\}_{i=1}^{I}.$ Together they form a Hamiltonian system with Hamiltonian given by $$ H_{\epsilon_l}\biggl(\{q_i\}_{i=1}^{I},\{p_i\}_{i=1}^{I}\biggr)=V_{\epsilon_l}\biggl(\{q_i\}_{i=1}^{I}\biggr) + \frac{p^T p}{2}, $$
where $p=(p_1,\ldots,p_{I})$ is the momenta vector. Thus, the momenta variables are Gaussian distributed random variables with mean zero and variance 1. The equations of motion for this Hamiltonian system are given by Hamilton's equations 
$$\frac{dq_i}{d\tau}=\frac{\partial H_{\epsilon_l}}{\partial p_i} \;\; \text{and} \;\; \frac{dp_i}{d\tau}=-\frac{\partial H_{\epsilon_l}}{\partial q_i} \;\; \text{for} \;\; i=1,\ldots,I.$$ 
Note that the Hamiltonian depends also on the initial condition for each sample $Y^n_0$ and we have written an equation for the evolution of $q_1=Y^n_0$ as well its associated momentum $p_1.$ Since the $Y^n_0$ are fixed by the resampling procedure we do not evolve them. However, note that the Brownian increment $\Delta B_0$ needs to be evolved because it affects the evolution of $Y_{\Delta t}.$ 

 HMC proceeds by assigning initial conditions to the momenta variables (through sampling  from $\exp(- \frac{p^T p}{2})$), evolving the Hamiltonian system in fictitious time $\tau$ for a given number of steps of size $\delta \tau$ and then using the solution of the system to perform a Metropolis accept/reject step (more details in \cite{liu}). After the Metropolis step, the momenta values are discarded. The most popular method for solving the Hamiltonian system, which is the one we also used, is the Verlet leapfrog scheme. In our numerical implementation, we did not attempt to optimize the performance of the HMC algorithm. For the sampling at each level of the drift relaxation process we used $10$ Metropolis accept/reject steps and $1$ HMC step of size $\delta \tau = 10^{-2}$ to construct a trial path. A detailed study of the drift relaxation/HMC algorithm for conditional path sampling problems outside of the context of particle filtering will be presented in a future publication.

\begin{figure}
\centering
\epsfig{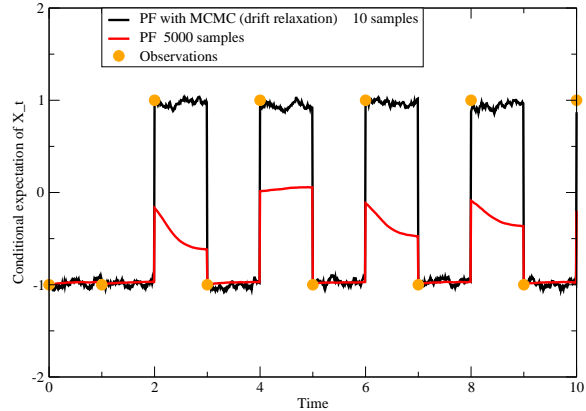}
\caption{Comparison of the conditional expectation of $X_t$ as computed by the generic particle filter and the particle filter with MCMC step. Different drift for each sample.}
\label{plot_comparison_varying}
\end{figure}

\begin{figure}
\centering
\epsfig{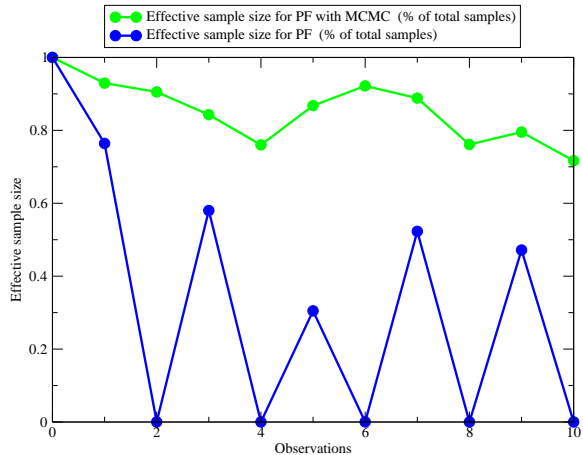}
\caption{Comparison of the effective sample size as computed by the generic particle filter and the particle filter with MCMC step. A different drift is used for each sample.}
\label{plot_comparison_ess_varying}
\end{figure}

For the chosen values of the parameters for the drift relaxation and HMC steps, the particle filter with MCMC step is about $500$ times more expensive per sample (particle) than the generic particle filter. However, we show that this increase in cost per sample is worthwhile. Figure \ref{plot_comparison_varying} compares the performance of the particle filter with MCMC step with $10$ samples and the generic particle filter with $5000$ samples. It is obvious that the particle filter with MCMC step follows accurately all the transitions between the two minima of the double-well. On the other hand, the generic particle filter captures accurately only every other observation. It fails to performs the transitions between the two minima of the double-well. Of course, since the particle filter with MCMC step uses only $10$ samples the conditional expectation estimate of the hidden signal is not as smooth as the estimate of the generic particle filter which uses $5000$ samples. However, the particle filter with MCMC step allows a much better resolution of the conditional density (conditioned on the observations). This can be seen by computing the effective sample size for the two filters. Figure \ref{plot_comparison_ess_varying} shows the effective sample size for both filters. Because of the different number of samples used in the two filters we have plotted the effective sample size for each filter as a percentage of the respective sample size. We see that the particle filter with MCMC step has overall a much better sample size than the generic particle filter. The generic particle filter has a wildly fluctuating effective sample size. In particular, since the generic particle filter misses every other observation, the corresponding effective sample size dips down to $1$ sample for the missed observations. Note from Figure \ref{plot_comparison_varying}  that even this one sample which dominates the observation weight does not come close to the observation. For the observations that the generic particle filter does capture, its effective sample size is still lower than the effective sample size of the particle filter with MCMC step. 

\section{Discussion}\label{discussion}
   
We have presented an algorithm for conditional path sampling of SDEs. The proposed algorithm is based on drift relaxation which allows to sample conditional paths from a modified drift equation. The conditional paths of the modified drift equation are then morphed into conditional paths of the original equation. We have called this process of gradually enforcing the drift of the original equation drift relaxation. The algorithm has been used to create a modified particle filter for SDEs. We have shown that the modified particle filter's performance is significantly better than the performance of a generic particle filter. 

In the current work, we have examined the application of drift relaxation to the filtering problem of diffusion in a double-well potential which is a standard example in the filtering literature. The same algorithm can be applied to the problem of tracking a single target. A problem of great practical interest is that of tracking not only one but multiple moving targets \cite{mahler,godsill,godsill2,vo}. The multi-target tracking problem is much more difficult than the single-target problem due to the combinatorial explosion of the number of possible target-observation association arrangements. In this context, the accurate tracking of each target becomes crucial. Suppose that only one of the targets is of interest and the rest act as decoys \cite{mahler2}. The inability to track each potential target accurately can lead to ambiguity about the targets' movement if the observations for different targets are close. We have already applied the drift relaxation modified particle filter to multi-target tracking problems with very encouraging results which will appear elsewhere \cite{MS10}.

\section*{Acknowledgements} 
I am grateful to Profs. V. Maroulas and J. Weare for many helpful discussions and comments. I would also like to thank for its hospitality the Institute for Mathematics and its Applications (IMA) at the University of Minnesota where the current work was completed.

\end{document}